\documentclass{elsart}
\usepackage{mathrsfs}
\usepackage{epsfig,latexsym,amsfonts,amsmath,amssymb}

\begin{document}
\begin{frontmatter}

\newtheorem{theorem}[subsection]{Theorem}
\newtheorem{lemma}[subsection]{Lemma}
\newtheorem{conjecture}[subsection]{Conjecture}
\newtheorem{proposition}[subsection]{Proposition}
\newtheorem{definition}[subsection]{Definition}
\newtheorem{corollary}[subsection]{Corollary}
\newtheorem{example}[subsection]{Example}
\newtheorem{remark}[subsection]{Remark}
\newtheorem{de}[subsection]{Definition}
\newtheorem{problem}[subsection]{Problem}
\newtheorem{question}[subsection]{Question}
\newtheorem{condition}[subsection]{Condition}

\renewcommand{\theequation}{\arabic{section}. \arabic{equation}}
\renewcommand{\thesection}{\arabic{section}}
\renewcommand{\thethm}{\arabic{section}.\arabic{thm}}
\renewcommand{\thefootnote}{\fnsymbol{footnote}}

\newcommand{\df}{\stackrel{\mbox{\rm def }}{=}}
\newcommand{\di}{\displaystyle}
\newcommand{\fix}{\mbox{Stab}}
\newcommand{\bl}[2]{{\left\langle #1 \:\:  \vrule \:\:  #2 \right\rangle}}
\newcommand{\vanish}[1]{}
\newcommand{\qbin}[2]{\left[\!\!\begin{array}{c}#1 \\ #2\end{array}\!\!\right]}

\title{Nontrivial independent sets of bipartite graphs and cross-intersecting families}

\author[wj]{Jun Wang}\ead{jwang@shnu.edu.cn}\ and
\author[wj,zh]{Huajun  Zhang\thanksref{Y}}\ead{huajunzhang@zjnu.cn}
\address[wj]{ Department of Mathematics, Shanghai Normal University,
Shanghai 200234, P.R. China}
\address[zh]{ Department of Mathematics,
 Zhejiang Normal University, Jinhua 321004, P.R. China}

\thanks[Y]{Corresponding author.  }

\baselineskip 20pt

\date{}
\maketitle
 \vspace{4mm}

\begin{abstract}
Let $G(X,Y)$ be a connected, non-complete bipartite graph with
$|X|\leq |Y|$. An independent set $A$ of  $G(X,Y)$ is said to be
trivial if $A\subseteq X$ or $A\subseteq Y$. Otherwise, $A$ is
nontrivial. By $\alpha(X,Y)$ we denote the size of maximal-sized
nontrivial independent sets of $G(X,Y)$. We prove that if the
automorphism group of $G(X,Y)$ is  transitive on $X$ and $Y$, then
$\alpha(X,Y)=|Y|-d(X)+1$, where $d(X)$ is the common degree of
vertices in $X$. We also give  the structures of maximal-sized
nontrivial independent sets of $G(X,Y)$.  As applications of this
result, we give the upper bound of sizes of two
cross-$t$-intersecting families of finite sets, finite vector spaces
and permutations.

\begin{keyword} intersecting family, cross-intersecting family,
symmetric system, Erd\H{o}s-Ko-Rado theorem
 \\[7pt]
\end{keyword}
\end{abstract}
\end{frontmatter}

\newcommand{\lr}[1]{\langle #1\rangle}
\newcommand{\qchoose}[2]{{ #1   \atopwithdelims[]  #2 }}

\

\parindent 17pt
\baselineskip 17pt
\section{Introduction}

Let $X$ be a finite  set and, for $0\leq k \leq |X|$, let
$\binom{X}{k}$ denote the family of all $k$-subsets of $X$, and let
$S_X$ and $A_X$  denote the symmetric group and alternative group on
$X$, respectively. In particular, for positive integer $n$, let
$[n]$ denote the set $\{1,2,\ldots,n\}$, $[k,n]=\{k+1,\ldots,n\}$
for $k\leq n$, $\overline  A=[n]\backslash A$ for $A\subseteq [n]$,
and abbreviate the symmetric group and alternative group on $[n]$ as
$S_n$ and $A_n$, respectively.

A family $\mathcal{A}$ of sets is said to be $t$-intersecting  if
$|A\cap B|\geq t$ holds for all $A,B\in\mathcal{A}$. Usually,
$\mathcal{A}$ is called intersecting if $t=1$. The celebrated
Erd\H{o}s--Ko--Rado theorem \cite{EKR},  says that if $\mathcal{A}$ is a
$t$-intersecting family in $\binom{[n]}{k}$, then
\[|\mathcal A|\leq \binom{n-t}{k-t}\]
for $n\geq n_0(k,t)$.
 The smallest $n_0(k,t)=
(k-t+1)(t+1)$ was determined by Frankl \cite{pf} for $t\geq 15$ and
subsequently determined by Wilson \cite{wilton} for all $t$.

The Erd\H{o}s-Ko-Rado theorem has many  generalizations, analogs and
variations. First, the notion of intersection is generalized to
t-intersection, and finite sets are analogous to finite vector
spaces, permutations and other mathematical objects.  Second,
intersecting families are generalized to cross-intersecting
families: $\mathcal A_1,\mathcal A_2, \ldots, \mathcal A_m$ are said
to be cross-$t$-intersecting if $|A\cap B|\geq t$ for all
$A\in\mathcal A_i$ and $B\in\mathcal A_j$,  $i\neq j$. Some typical
but far from exhaustive results are listed as follows.

Let $\mathbb F_q$ be a finite field of order $q$, $V=V_n(q)$ an
$n$-dimensional vector space over $\mathbb F_q$ , and
$\qchoose{V}{k}$ the set of all $k$-dimensional subspaces (or
$k$-subspace, for short) of $V$. Then the cardinality of
$\qchoose{V}{k}$ equals $\qchoose
nk_q=\prod_{i=0}^{k-1}\frac{q^{n-i}-1}{q^{k-i}-1}$. For brevity, we
write 
 $\qchoose nk$ rather than
 $\qchoose nk_q$.
 A subset
$\mathcal A$ of $\qchoose{V}{k}$ is said to be a $t$-intersecting
family if $\dim(A\cap B)\geq t$ for any $A,B\in \mathcal A$. The
Erd\H{o}s-Ko-Rado theorem for finite vector spaces says that if
$\mathcal{A}$ is a $t$-intersecting family in$\qchoose{V}{k}$, then
\[|\mathcal A|\leq \max\left\{\qchoose
{n-t}{k-t},\qchoose{2k-t}{k}\right\}\] for $n\geq 2k-t$. This
theorem was  first established by Hsieh \cite{Hsieh} for $t=1$ and
$k<n/2$,  then by Greene and  Kleitman \cite{Greene-Kleit} for $t=1$
and $k|n$, and finally  by Frankl and Wilson \cite{Frankl-Wilson}
for the general case.

 A subset $A$ of $S_n$ is said to be a
$t$-intersecting family if any two permutations in $A$ agree in at
least $t$ points, i.e. for any $\sigma,\tau\in A$, $|\{i\in
[n]:\sigma(i)=\tau(i)\}|\geq t$.  Deza and Frankl \cite{DF} showed
that an intersecting family in $S_n$ has size at most $(n-1)!$ and
conjectured that for $t$ fixed, and $n$ sufficiently large depending
on $t$, a $t$-intersecting family in $S_n$ has size at most
$(n-t)!$. Cameron and Ku \cite{Cameron} proved an intersecting
family of size  $(n-1)!$ is a coset of the stabilizer of a point. A
few alternative proofs of Cameron and Ku's result are given in
\cite{Larose}, \cite{godsil} and \cite{wz}. Ku and Leader \cite{Ku}
also generalized this result to partial permutations (see also
\cite{lw}).  Ellis, Friedgut and Pilpel  \cite{Ellis} proved Deza
and Frankl's conjecture on $t$-intersecting family in $S_n$.

 Hilton \cite{hilton} investigated the cross-intersecting families in
$\binom{[n]}{k}$: Let
$\mathcal{A}_1,\mathcal{A}_2,\ldots,\mathcal{A}_m$ be
cross-intersecting families in $\binom{[n]}{k}$ with
$\mathcal{A}_1\neq\emptyset$. If  $k\leq n/2$, then
\begin{equation}\label{hilton}
 \sum_{i=1}^m|\mathcal{A}_i|\leq\left\{
                                        \begin{array}{cl}
                                          \binom{n}{k}, & \hbox{if $m\leq \frac{n}{k}$;} \\
                                          m\binom{n-1}{k-1}, & \hbox{if $m\geq \frac{n}{k}$.}
                                        \end{array}
                                      \right.
\end{equation}
He also determined the structures of $\mathcal{A}_i$'s when equality
holds. Borg \cite{pborg1} gives a simple proof of this theorem, and
generalizes it to   labeled sets \cite{borg}, signed sets \cite{borg6} and permutations
\cite{pborg}. We generalized this theorem to general symmetric
systems \cite{wzhj}, which contain finite sets, finite vector spaces
and permutations, etc.

 Hilton and Milner \cite{hm} and Frankl and Tokushige \cite{FT} also investigated
  the sizes of two cross-intersecting families:
If $\mathcal{A}\subset \binom{[n]}{a}$ and
$\mathcal{B}\subset\binom{[n]}{b}$ are  cross-intersecting families
with $n\geq a+b$, $a\leq b$, then $|\mathcal{A}|+|\mathcal{B}|\leq
\binom{n}{b}-\binom{n-a}{b}+1$.

This theorem actually gives a upper bound of the sizes of nontrivial
independent sets in a bipartite graph.

Let $G$ be a simple graph with vertex set $V(G)$ and edge set
$E(G)$. For $v\in V(G)$, define $N_G(v)=\{u\in V(G):uv\in E(G)\}$
and $N_G(A)=\cup_{v\in A}N_G(v)$ for $A\subseteq V(G)$. If there is
no possibility of confusion, we abbreviate $N_G(A)$ as $N(A)$. A
subset $A$ of $V(G)$ is an independent set of $G$ if $A\cap
N(A)=\emptyset$. A graph $G$ is bipartite if $V(G)$ can be
partitioned into two subsets $X$ and $Y$ so that every edge has one
end in $X$ and one end in $Y$. In this case, we denote the bipartite
graph by $G(X,Y)$.
 An independent set $A$ of $G(X,Y)$ is said to be trivial if $A\subseteq X$ or $A\subseteq Y$.
 In any other case, $A$ is nontrivial.
 If every vertex
in $X$ is adjacent to every vertex in $Y$, then $G(X,Y)$ is called a
complete bipartite graph. Clearly, a complete bipartite graph has
only trivial independent sets. A bipartite graph $G(X,Y)$ is said to
be {\em part-transitive} if there is a group $\Gamma$ transitively
acting on $X$ and $Y$, respectively, and preserving the adjacency
relation of the graph. Clearly, if $G(X,Y)$ is part-transitive, then
every vertex of $X$ ($Y$) has the same degree, written as $d(X)$
($d(Y)$). By $\alpha(X,Y)$ and $I(X,Y)$ we denote the size and the
set of maximal-sized nontrivial independent sets of $G(X,Y)$,
respectively.


This paper contributes to $\alpha(X,Y)$ and $I(X,Y)$ for
part-transitive bipartite graphs $G(X,Y)$. To do this we make a
simple observation as follows.

Let $G(X,Y)$ be a non-complete bipartite graph and let $A\cup B$ be
a nontrivial independent set of $G(X,Y)$, where $A\subset X$ and
$B\subset Y$. Then $A\subseteq X\backslash  N(B)$ and $B\subseteq
Y\backslash N(A)$, which implies that
\[|A|+|B|\leq \max\{|A|+|Y|-|N(A)|, |B|+|X|-|N(B)|\}.\]
From this one sees that
\begin{equation}\label{maxi}
\alpha(X,Y)=\max\{|Y|-\epsilon(X), |X|- \epsilon(Y)\},
\end{equation}
where $$ \epsilon(X)=\min\{|N(A)|-|A|: \mbox{$A\subset X$ and
$N(A)\neq Y$}\}$$ and
$$\epsilon(Y)=\min\{|N(B)|-|B|: \mbox{$B\subset Y$ and $N(B)\neq X$}\}.$$

A subset $A$ of $X$  is called a \emph{fragment} in $X$ if $N(A)\neq
Y$  and $|N(A)|-|A|=\epsilon(X)$.  By $\mathcal F(X)$ we denote the
set of all fragments contained in $X$. $\mathcal F(Y)$ is defined in
a similar way and write  $\mathcal F(X,Y)=\mathcal F(X)\cup \mathcal
F(Y)$. An element  $A\in \mathcal F(X,Y)$ is also called a
$k$-fragment if $|A|=k$. As we shall see (Lemma \ref{lad}) that
$|Y|-\epsilon(X)=|X|- \epsilon(Y)$. Therefore,  in order to address
our problems it suffices to determine $\mathcal F(X)$ or $\mathcal
F(Y)$.

Let $X$ be a finite set, and $\Gamma$ a group transitively acting on
$X$. We say the action of $\Gamma$ on $X$ is \emph{primitive}, or
$\Gamma$ is primitive on $X$,  if $\Gamma$ preserves no nontrivial
partition of $X$. In any other case, the action of $\Gamma$ is
\emph{imprimitive}. It is easy to see that  if the action of
$\Gamma$ on $X$ is transitive and imprimitive, then there is a
subset $B$ of $X$ such that $1<|B|<|X|$ and $\gamma(B)\cap B=B$ or
$\emptyset$ for every $\gamma\in G$. In this case, $B$ is called an
\emph{imprimitive set}  in $X$. It is well known that the action of
$\Gamma$ is primitive if and only if for each $a\in X$, the
stabilizer of $a$, written as $\Gamma_a$ defined to be the set
$\{\gamma\in \Gamma:\gamma(a)=a\}$, is a maximal subgroup of
$\Gamma$ (cf. \cite[Theorem 1.12]{BA}). Furthermore, a subset $B$ of
$X$ is said to be \emph{semi-imprimitive}  if $1< |B|<|X|$ and
$|\gamma(B)\cap B|=0,1$ or $|B|$ for each $\gamma\in \Gamma$.
Clearly,  every $2$-subset of $X$ is semi-imprimitve.

The following are main results of this paper.

\begin{theorem}\label{Gd}
Let $G(X,Y)$ be a   non-complete bipartite graph with $|X|\leq |Y|$.
If $G(X,Y)$ is part-transitive  and every fragment in $X$ and $Y$ is
primitive under the action of  a group $\Gamma$. Then
$\alpha(X,Y)=|Y|-d(X)+1$. Moreover,
\begin{enumerate}
\item [\rm(i)] if $|X|<|Y|$, then each fragment in $X$ has size $1$;
\item[\rm(ii)] if $|X|=|Y|$, then
each fragment in $X$  has size $1$ or $|X|-d(X)$ unless there is a
semi-imprimitive fragment in  $X$ or $Y$.
\end{enumerate}
\end{theorem}

As consequences of this theorem we give the upper bounds of sizes of
two cross-$t$-intersecting families of finite sets, finite vector
spaces and symmetric groups.

\begin{theorem}\label{subsets}
Let $n,a,b,t$ be positive integers  with $n\geq 4$, $a,b\geq 2$, $t<
\min\{a,b\}$, $a+b<n+t$, $(n,t)\neq (a+b,1)$ and $\binom{n}{a}\leq
\binom{n}{b}$. If $\mathcal{A}\subset\binom{[n]}{a}$ and
$\mathcal{B}\subset\binom{[n]}{b}$ are cross-t-intersecting, then
\begin{equation}\label{cross1}
|\mathcal{A}|+|\mathcal{B}|\leq
\binom{n}{b}-\sum_{i=0}^{t-1}\binom{a}{i}\binom{n-a}{b-i}+1.
\end{equation}
 Moreover,
\begin{enumerate}
  \item [\rm(i)] when $\binom na<\binom nb$, equality holds if and only if $\mathcal{A}=\{A\}$ and
  $\mathcal B=\binom{[n]}{b}\backslash N(A)$ for any $A\in\binom{[n]}{a}$;
  \item [\rm(ii)] when $\binom na=\binom nb$,  equality holds if and only if  either $\mathcal{A}=\{A\}$
  and $\mathcal B=\binom{[n]}{b}\backslash N(A)$ for any $A\in\binom{[n]}{a}$, or  $\mathcal{B}=\{B\}$ and $\mathcal A=\binom{[n]}{a}
  \backslash N(B)$ for any   $B\in\binom{[n]}{b}$,
  or
  $\{a,b,t\}=\{2,2,1\}$ and $\mathcal{A}=\mathcal{B}=\{C\in \binom{[n]}{2}: i\in C\}$ for some $i\in [n]$, or
   $\{a,b,t\}=\{n-2,n-2,n-3\}$ and $\mathcal{A}=\mathcal{B}=\binom{A}{n-2}$ for some $A\in \binom{[n]}{n-1}$.
  \end{enumerate}
\end{theorem}

\begin{theorem}\label{subspaces} Let $V$  be an $n$-dimensional vector space over
the   field of order $q$ and let $n,a,b,t$ be positive integers with
$n\geq 4$, $a,b\geq 2$, $t< \min\{a,b\}$, $a+b<n+t$,  and
$\qchoose{n}{a}\leq \qchoose{n}{b}$.  If  $\mathcal{A}\subset
\qchoose{V}{a}$ and $\mathcal{B}\subset\qchoose{V}{b}$  are
cross-$t$-intersecting,
then
\begin{equation}\label{cross2}
|\mathcal{A}|+|\mathcal{B}|\leq
\qchoose{n}{b}-\sum_{i=0}^{t-1}q^{(a-i)(b-i)}\qchoose{a}{i}\qchoose{n-a}{b-i}+1.
\end{equation}
 Moreover, equality holds if and only if
 $\mathcal{A}=\{A\}$ and $\mathcal{B}= \qchoose{V}{b}\backslash N(A)$ where $A\in \qchoose{V}{a}$, or
 $\mathcal{A}=\qchoose{V}{b}\backslash N(B)$ and $\mathcal{B}=\{B\}$ where $B\in  \qchoose{V}{b}$, subject to $\qchoose{n}{a}=
\qchoose{n}{b}$.
\end{theorem}

\begin{theorem}\label{group} Let $n$ and $t$ be positive integers  with  $n\geq 4$ and $t\leq n-2$.
If  $\mathcal{A}$ and $\mathcal{B}$  are cross-$t$-intersecting
families in $S_n$, then
\begin{equation}\label{cross3}
|\mathcal{A}|+|\mathcal{B}|\leq n!-\sum_{i=0}^{t-1}{n\choose
i}D_{n-i}+1,
\end{equation}
where $D_{n-i}$ is the number of derangements in $S_{n-i}$.
 Moreover, equality holds if and only if
 $\{\mathcal{A},\mathcal{B}\}=\left\{\{\sigma\}, S_n\backslash N(\sigma)\right\}$ where $\sigma\in S_n$.
\end{theorem}

We shall prove Theorem \ref{Gd} in the next section, Theorem
\ref{subsets} in Section 3, Theorem \ref{subspaces} in Section
4  and Theorem \ref{group} in Section
5.

\section{Proof of Theorem \ref{Gd}}
Before to start the proof of Theorem \ref{Gd} we present two lemmas.

\begin{lemma}\label{lad}
Let $G(X,Y)$ be a non-complete bipartite graph. Then,
$|Y|-\epsilon(X)=|X|-\epsilon(Y)$, and
\begin{enumerate}
\item [\rm(i)] $A\in
\mathcal F(X)$ if and only if $Y\backslash N(A)\in \mathcal F(Y)$,
and $N(Y\backslash N(A))=X\backslash A$;
\item [\rm(ii)]$A\cap B$
and $A\cup B$ are both in $\mathcal F(X)$ if $A,B\in\mathcal F(X)$,
$A\cap B\neq \emptyset$ and $N(A\cup B)\neq Y$.

\end{enumerate}
\end{lemma}
\textbf{Proof.} Suppose $A\in \mathcal F(X)$ and put  $C=Y\backslash
N(A)$. Clearly, $N(C)\subseteq X\backslash A$. If $N(C)\neq
X\backslash A$, writing $A'=X\backslash N(C)$, then $A\subsetneq A'$
and $N(A')=N(A)$. So $|N(A')|-|A'|<|N(A)|-|A|=\epsilon(X)$, yielding
a contradiction. Hence $N(C)=X\backslash A$, and
$|N(C)|-|C|=(|X|-|A|)-(|Y|-|N(A)|)=\epsilon(X)-|Y|+|X|\geq\epsilon(Y)$.
Symmetrically, for $D\in \mathcal F(Y)$, putting $A=X\backslash
N(D)$, we have $N(A)=Y\backslash D$ and
$|N(A)|-|A|=(|Y|-|D|)-(|X|-|N(D)|)=\epsilon(Y)-|X|+|Y|\geq\epsilon(X)$.
We then obtain that $\epsilon(X)+|X|=\epsilon(Y)+|Y|$ and (i) holds.

Now, suppose that $A,B\in \mathcal F(X)$,  $A\cap B\neq \emptyset$
and $N(A\cup B)\neq Y$. Then $ |N(A\cup B)|-|A\cup
B|\geq\epsilon(X)$ and $|N(A\cap B)|-|A\cap B|\geq \epsilon(X)$.
Note that $N(A\cup B)=N(A)\cup N(B)$ and $N(A\cap B)\subseteq
N(A)\cap N(B)$. We have
  \begin{eqnarray*}
\epsilon(X)&\leq& |N(A\cup B)|-|A\cup
B|\\&=&|N(A)|+|N(B)|-|N(A)\cap N(B)|-|A|-|B|+|A\cap B|\\
&\leq& 2\epsilon(X)-(|N(A\cap B)|-|A\cap B|)\leq\epsilon(X),
\end{eqnarray*}
which implies that $|N(A\cup B)|-|A\cup B|=\epsilon(X)$ and
$|N(A\cap B)|-|A\cap B|=\epsilon(X)$, hence (ii) holds. \qed

From the first statement of this lemma it follows that there is a
one to one correspondence  $\phi: \mathcal F(X,Y)\mapsto \mathcal
F(X,Y)$, where
\[\phi(A)=\left\{\begin{array}{ll}
Y\backslash N(A) & \mbox{if $A\in \mathcal F(X)$},\\
X\backslash N(A)& \mbox{if $A\in \mathcal F(Y)$}.
\end{array}\right.\] Moreover, $\phi$ is an involution, i.e., $\phi^{-1}=\phi$, and
$|A|+|\phi(A)|=\alpha(X,Y)$.  A fragment is called {\em balanced} if
$|A|=|\phi(A)|$. Clearly, all balanced fragments have identical size
$\frac 12\alpha(X,Y)$.



\begin{lemma}\label{lad1} Let $G(X,Y)$ be a non-complete and part-transitive bipartite
 graph under the action of a group
$\Gamma$. Suppose that $A\in \mathcal F(X,Y)$ such that
$\emptyset\neq\gamma(A)\cap A\neq A$ for some $\gamma\in \Gamma$. If
$|A|\leq |\phi(A)|$, then $A\cup\gamma(A)$ and $A\cap \gamma(A)$ are
both in $\mathcal F(X,Y)$.
\end{lemma}
\textbf{Proof.} Without loss of generality, suppose $A\in\mathcal
F(X)$ and $|A|\leq |\phi(A)|=|Y\backslash N(A)|$.  Since
$|N(A)|=|A|+\epsilon(X)$ and $|N(A\cap \gamma(A))|\geq |A\cap
\gamma(A)|+\epsilon(X)$,
\begin{eqnarray*}
|N(A\cup \gamma(A))| &=&2|N(A)|-|N(A)\cap N(\gamma(A))|\\
 &\leq& 2 |N(A)|-|N(A\cap \gamma(A))|\\
\nonumber &\leq&
|N(A)|+|A|+\epsilon(X)-(|A\cap \gamma(A)|+\epsilon(X))\\
&=&|N(A)|+|A\backslash \gamma(A)|<|N(A)|+|Y\backslash
N(A)|=|Y|.\end{eqnarray*} Then, by Lemma \ref{lad} (ii),  $A\cap
\gamma(A)$ and $A\cup \gamma(A)$ are both in $\mathcal F(X)$.\qed

From the above lemma  we have that if every element of  $X$ ($Y$) is
primitive and there is an $A\in\mathcal F(X)$ ($\mathcal F(Y)$) with
$|A|\leq |\phi(A)|$, then $\mathcal F(X)$ ($\mathcal F(Y)$) contains
a singleton. In particular, when $|X|=|Y|$  there are always two
kinds of fragments in $X$: one is $\{a\}$ for $a\in X$, the other is
$X\backslash N(b)$ for $b\in Y$. The former is a minimal-sized
fragment, and the latter is maximal-sized one. We call the fragments of
this kinds {\it trivial}. All others are  {\it
nontrivial}. 

\textbf{Proof of Theorem \ref{Gd}.} From the above discussion we
have that  $\mathcal F(X)$ or $\mathcal F(Y)$ contains a singleton,
that is, $\alpha(X,Y)=\max\{|Y|-d(X)+1, |X|-d(Y)+1\}$. By counting
the  edges of  $G(X,Y)$  we have $d(X)|X|=d(Y)|Y|$, so
$d(X)=d(Y)|Y|/|X|\geq d(Y)$ because $|Y|\geq |X|$. Then
\[|Y|-|X|=d(X)|X|/d(Y)-|X|=(d(X)-d(Y))|X|/d(Y)\geq d(X)-d(Y).\]
Equality holds if and only if $d(X)=d(Y)$ hence $|X|=|Y|$ because
$|X|>d(Y)$.  This proves that $|X|-d(Y)+1\leq |Y|-d(X)+1$ and
equality holds if and only if $|X|=|Y|$. In any cases,
$\alpha(X,Y)=|Y|-d(X)+1$.

We complete the proof by two cases.

Case 1: $|X|<|Y|$.  In this case we have seen that $\mathcal F(X)$
contains singletons while $\mathcal F(Y)$ does not. Now, let $A$ be
a maximal-sized element of $\mathcal F(X)$ and write
$B=\phi(A)=Y\backslash N(A)$.
 Then $B$ is a minimal-sized element of $\mathcal F(Y)$ with  $|B|>1$ and $\phi(B)=A$.
  Suppose $|A|>1$.
Since  $A$ and $B$  are primitive,  there are  $\sigma, \gamma\in
\Gamma$ such that $\sigma(A)\cap A\neq \emptyset$, $\sigma(A)\neq
A$, $\gamma(B)\cap B\neq \emptyset$ and $\gamma(B)\neq B$. From this
and Lemma \ref{lad1} it follows that  if $A|\leq |B|=|\phi(A)|$,
then $\sigma(A)\cup A\in \mathcal F(X)$, contradicting the
maximality of $|A|$; if $|B|\leq |A|=|\phi(B)|$, then $\gamma(B)\cap
B\in \mathcal F(Y)$, contradicting the minimality of $|B|$. This
proves that  $|A|=1$ for every $A\in \mathcal F(X)$.

Case 2: $|X|=|Y|$. In this case,  if there is a nontrivial fragment
in $X$ or in $Y$, let $A$ be a minimal-sized one. Then $1<|A|\leq
|\phi(A)|$. From Lemma \ref{lad1} it follows that
 for every $\gamma\in \Gamma$,  $\gamma(A)\cap A$ is a fragment whenever
$\gamma(A)\cap A\neq \emptyset$. Then, the minimality of $|A|$
implies that $|\gamma(A)\cap A |=0$, 1 or $|A|$, for every
$\gamma\in \Gamma$, i.e.,  $A$ is  semi-imprimitve.

The proof is complete.  \qed

For applications of the theorem, we make further discussions on the
fragments in the rest of this section.

\begin{figure}[hh]\label{f1}
 \centering
     \includegraphics[width=1.5 in]{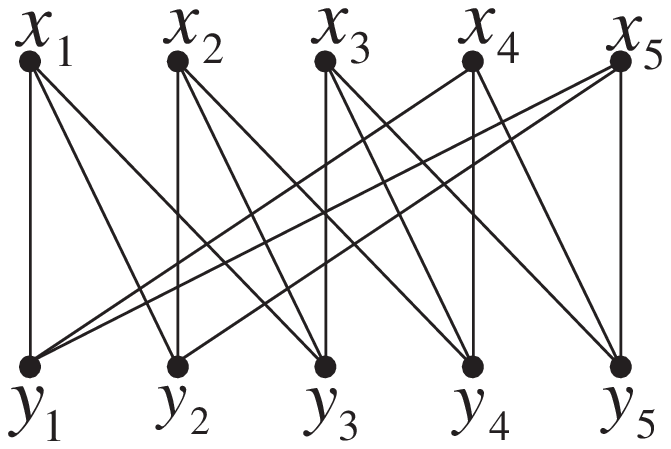}
  \caption{}
\end{figure}

 Note that most bipartite graphs concerning here
have only trivial fragments, but there are actually bipartite
graphs, which have sufficiently large nontrivial fragments. For
example, let $n$ and $r$ are fixed positive integer with $r<n$,
$X=\{x_1,x_2,\ldots,x_n\}$ and $Y=\{y_1,y_2,\ldots,y_n\}$. Define
$x_iy_j$ to be  an edge of $G(X,Y)$ if and only if $j\in\{i,
i+1,\ldots, i+r-1\} \pmod n$ (see Fig. \ref{f1} for $n=5$ and
$r=3$).
It is easy to verify that $\{x_i,x_{i+1},\ldots,x_{i+j}\}$, where
$1\leq j\leq n-r-1$ and the subscripts are computed modulo $n$,  is
a fragment in  $X$.

However, as we shall see, whether or not  a bipartite graph has
sufficiently large fragments depends  if it has a 2-fragment.
\begin{proposition}\label{2frag}
Let $G(X,Y)$ be a   non-complete  bipartite graph with $|X|=|Y|$ and
$\epsilon(X)=d(X)-1$, and let $\Gamma$ be a group  part-transitively
acting on  $G(X,Y)$. If there is a 2-fragment in $X$, then either
\begin{enumerate}
  \item [\rm(i)]  there is an imprimitive set $A\subset X$  with
$|N(A)|-|A|=d(X)-1$, or
  \item [\rm(ii)] there is a subset $A\subseteq X$, where $A=X$ or $A$ is an
imprimitive set under the action of $\Gamma$ with $|A|>2$, such that
the quotient group $\Gamma_A/(\cap_{a\in A}\Gamma_a)$ is isomorphic
to a subgroup of the dihedral group $D_{|A|}$, where
$\Gamma_A=\{\sigma\in\Gamma:\sigma(A)=A\}$.
      \end{enumerate}
\end{proposition}
\textbf{Proof.} By definition we have that for any  $x,y\in X$,
$\{x,y\}$  is a 2-fragment if and only if $|N(x)\cap N(y)|=d(X)-1$.
We now define a simple graph $H(X)$, whose vertex set $V(H)$ is $X$,
and whose edge set $E(H)$ consists of all pairs $xy$'s such that
$\{x,y\}$ is a fragment in $X$. Then, each element of $\Gamma$
induces an automorphism of $H(X)$. So $H(X)$ is vertex-transitive.
As usual, the valency of $H(X)$ is denoted by $d(H)$.

 Let $H'$ be a connected
component of $H(X)$ and let $A$ be the vertex set of $H'$.  Then
$|A|\geq 2$. If $|A|=2$, then $A$ is clearly an imprimitive set in
$A$ with $|N(A)|-|A|=d(X)-1$.  Suppose  $|A|>2$ and let  $xyz$ be a
path in $H'$ for distinct $x,y,z\in A$. Set $N(x)=B\cup \{a\}$ and
$N(y)=B\cup \{b\}$, where $B=N(x)\cap N(y)$. Since $yz\in E(H')$,
$N(z)=(N(y)\setminus\{c\})\cup\{d\}$ for some $c\in N(y)$ and $d\in
N(z)\setminus N(y)$. If  $d\in N(x)\cup N(y)$, then
$N(\{x,y,z\})=N(\{x,y\})$, contradicting that $\{x,y\}$ is a
fragment. Therefore, $d\not\in N(x)\cup N(y)$. So
$$N(x)\cap N(z)=\left\{
                      \begin{array}{cl}
                        B, & \hbox{if $c=b$,} \\
                        B\setminus\{c\}, & \hbox{if $c\neq b$.}
                      \end{array}
                    \right.
$$
From this it follows immediately the following.

{\em Claim}: $N(y)\subset N(x)\cup N(z)$ if  the induced subgraph
$H'[\{x,y,z\}]$ is a path, and $|N(x)\cap N(y)\cap N(z)|=d(X)-1$ if
$H'[\{x,y,z\}]$  is a cycle.

 If $H'$ is a complete graph, then from the above claim  it follows that $|\cap_{x\in A}N(x)|=d(X)-1$,
 so $|N(A)|-|A|=d(X)-1$. Since $d(X)\geq 2$, we have $|A|<|X|$,  hence $A$ is an imprimitive set in  $X$ with
 $|N(A)|-|A|=d(X)-1$.

If  $H'$ is not complete, then there are more than three elements of
$A$, say $x_1,x_2,\ldots,x_m$,  such that the induced subgraph
 $H'[\{x_1,\ldots,x_m\}]$ is a cycle, written as $x_1x_2\cdots x_mx_1$.
By definition we see that $|N(\{x_1,x_{2},\ldots,x_{s}\})|\leq
d(X)-1+s$ for $1\leq s\leq m-1$, and equality holds if
$N(\{x_1,x_{2},\ldots,x_{s}\})\neq Y$, that is,
$\{x_1,x_{2},\ldots,x_{s}\}$ is a fragment. Now, if
$N(\{x_1,x_{2},\ldots,x_{m-1}\})\neq Y$, then, by the above claim,
$N(\{x_1,x_{2},\ldots,x_{m}\})=N(\{x_1,x_{2},\ldots,x_{m-1}\})\neq
Y$, which yields a contradiction since
$|N(\{x_1,x_{2},\ldots,x_{m}\})|-m<d(X)-1$. Therefore,
$N(\{x_1,x_{2},\ldots,x_{m-1}\})=Y$. Assume that $t$ is the least
index such that $N(\{x_1,x_{2},\ldots,x_{t}\})=Y$, where $2<t\leq
m-1$. This means that every path of length less than $t$ on this
cycle is a fragment. In this case, if $d(H)>2$, then there is an
$x\in A\backslash\{x_1,\ldots,x_m\}$ such that $xx_{t-1}$ is an edge
of $H'$. Setting $a\in N(x)\backslash N(x_{t-1})$, we have that
$a\in N(x_i)$ for some $i\in [t]\backslash\{t-1\}$. Then
$N(\{x_i,\ldots,x_{t-1}\})=N(\{x_i,\ldots,x_{t-1},x\})$ if $i<t-1$,
or $N(\{x_{t-1}, x_t\})=N(\{x_{t-1},x_t,x\})$ if $i=t$.  Both the
cases contradict that $\{x_i,\ldots,x_{t-1}\}$ and
$\{x_{t-1},x_{t}\}$ are fragments. This proves that $H'$ is a cycle,
and hence (ii) holds. \qed

 \begin{proposition}\label{3}
Let $G(X,Y)$ be as in Proposition \ref{2frag}. If there are no
2-fragments in $\mathcal F(X,Y)$, then every nontrivial fragment
$A\in \mathcal F(X)$  (if it exists) is balanced,  and for each
$a\in A$, there is a unique nontrivial  fragment $B$ such that
$A\cap B=\{a\}$.
\end{proposition}
\textbf{Proof.}  Let $A$ be a minimal-sized nontrivial fragment in
$X$ and $Y$. Then, $\phi(A)$ is a maximal-sized fragment in $X$ and
$Y$. Without loss of generality, suppose that $A\in \mathcal{F}(X)$.
Then $\phi(A)=Y\backslash N(A)$ and  $|Y\backslash N(A)|\geq |A|$.
We now prove that the equality holds, i.e., $A$ is balanced.
Suppose, to the contrary, that $|Y\backslash N(A)|>|A|$. Set
$\mathcal{A}=\{\sigma(A): \sigma\in \Gamma\}$.
 As we have mentioned, $A$ is semi-imprimitive, so $|B\cap C|=1$ or $0$
 for all
distinct  $B, C\in \mathcal{A}$.   We now define a graph
$H(\mathcal{A})$, whose vertex set is $\mathcal{A}$, and whose edge
set consists of all pairs $BC$'s such that $|B\cap C|=1$ for
$A,B\in\mathcal A$. Clearly, $H(\mathcal{A})$ is vertex-transitive.
Since $A$ is primitive, $H(\mathcal A)$ is not an empty graph.
Suppose  that $A\cap B=\{b\}$ for some $B\in\mathcal A$ and $b\in
A$.
 Then,  for each $a\in A$, the part-transitivity of $G(X,Y)$ implies that there is a $\sigma\in \Gamma$ with
 $\sigma(b)=a$, hence $\sigma(A)\cap\sigma(B)=\{a\}$. From this it follows that
 the valency of $H(\mathcal{A})$, denoted by $d(H)$, is at
least $|A|>2$. Hence $H(\mathcal{A})$ contains a cycle. Let
$AA_1\ldots A_sA$ be one of minimum length. Then the induced
subgraph $H[\{A, A_1,\ldots,A_i\}]$ is a path from $A$ to $A_i$ for
$i=1,2,\ldots, s-1$. By Lemma \ref{lad}, if $N(A\cup
A_1\cup\cdots\cup A_i)\neq Y$, then both $A\cup A_1\cup\cdots\cup
A_i$ and $Y\backslash N(A\cup A_1\cup\cdots\cup A_i)$ are fragments.
Furthermore, if $|Y|-|N(A\cup A_1\cup\cdots\cup A_{i})|>1$, then the
minimality of $A$ implies $|Y|-|N(A\cup A_1\cup\cdots\cup
A_{i})|\geq |A|$, hence
\begin{eqnarray*} |N(A\cup A_1\cup\cdots\cup A_{i+1})| &\leq&
|N(A\cup A_1\cup\cdots\cup A_{i})|+|N(A_{i+1})|\\ && -|N((A\cup
A_1\cup\cdots\cup A_{i})\cap A_{i+1})|\\
&=& |N(A\cup A_1\cup\cdots\cup A_{i})|+|A|-1\leq |Y|-1,
\end{eqnarray*}
i.e., $A\cup A_1\cup\cdots\cup A_{i+1}$ is also a fragment. Now, if
$|Y|-|N(A\cup A_1\cup\cdots\cup A_{s-1})|>1$, then, by Lemma
\ref{lad}, $A_s\cap (A\cup\cdots\cup A_{s-1})$ is a fragment.
However, it is clear that $|A_s\cap (A\cup\cdots\cup A_{s-1})|=2$,
yielding  a contraction. Therefore, there is a unique index $k$ with
$2\leq k\leq s-1$ such that  $|Y|-|N(A\cup A_1\cup\cdots\cup
A_{k})|=1$, that is, $A\cup A_1\cup\cdots\cup A_{k}$ is a
maximal-sized fragment. In this case, it is clear that
 $|Y|-|N(A\cup
A_1\cup\cdots\cup A_{k-1})|=|A|$. We now find a contradiction.


Set $A'=A\cup A_1\cup\cdots\cup A_{k-1}$. Then  for each $B\in
(N_H(A)\cup N_H(A_{k-1}))\setminus \{A_1,A_{k-2}\}$, the induced
subgraph  $H[\{A, A_1,\ldots, A_{k-1},B\}]$ is a path of length $k$
in $H(\mathcal A)$, so the above argument is available here.   We
thus obtain  at least $2(d(H)-1)$ many maximal-sized fragments in
$X$ containing $A'$.  On the other hand, for every maximal-sized
fragment $C\in \mathcal{F}(X)$ containing $A'$, we have that
$C=X\backslash N(b)$ for some $b\in Y\backslash N(A')$  since
 $|Y\backslash N(C)|=1$,   hence there
are at most $|Y\backslash N(A')|=|A|$ many maximal-sized fragments
in $\mathcal F(X)$  containing  $A'$, yielding  a contradiction
because $2(d(H)-1)\geq 2(|A|-1)>|A|$. This proves that $|Y\backslash
N(A)|=|A|$, i.e., $A$ is balanced. Assume $A=\{a_1,\ldots,a_d\}$
where $d>2$. As we have seen, for each $i$, there is a $\sigma_i\in
\Gamma$ such that $A\cap\sigma_i(A)=\{a_i\}$. Then
$A\cup\sigma_i(A)$ is a maximal-sized fragment containing $A$, and
the semi-imprimitivity and $d>2$ imply  $A\cup\sigma_i(A)\neq
A\cup\sigma_j(A)$ if $i\neq j$. Therefore, $A\cup\sigma_i(A)$,
$i=1,\ldots,d$, are the all maximal-sized fragments containing $A$.
This proves that for every $a\in A$, there is only one nontrivial
fragment $B$ with $A\cap B=\{a\}$.
 \qed

\section{Proof of Theorem \ref{subsets}}
With  the assumptions in the theorem, we put $\mathcal
{X}=\binom{[n]}{a}$ and $\mathcal {Y}=\binom{[n]}{b}$. The bipartite
graph $G(\mathcal{X},\mathcal{Y})$ is defined by the
cross-$t$-intersecting relation between  $\mathcal{X}$ and
$\mathcal{Y}$:  For $A\in\mathcal{X}$ and $B\in\mathcal{Y}$, $AB\in
E(G)$ if and only if $|A\cap B|<t$. It is easy to check that
$G(\mathcal{X},\mathcal{Y})$ is connected since $a+b<t+n$ and
$(n,t)\neq (a+b,1)$, and  $G(\mathcal{X},\mathcal{Y})$ is
non-complete since $t<\min\{a,b\}$. Clearly,
 $S_n$ transitively acts on $\mathcal {X}$ and $\mathcal {Y}$, respectively, in a natural way, and preserves the
 cross-$t$-intersecting relation. Therefore,  $d(\mathcal X)=|N(A)|$
 for each $A\in \mathcal X$, and $d(\mathcal Y)=|N(B)|$
 for each $B\in \mathcal Y$.  It is easy to see that
for each $A\in \mathcal {X}$,
 $$N(A)=\{B\in \binom{[n]}{b}: |A\cap B|<t\}=\bigcup_{0\leq
i\leq t-1}\{B\in \binom{[n]}{b}: |A\cap B|=i\},$$
 hence $|N(A)|=\sum_{i=0}^{t-1}\binom{a}{i}\binom{n-a}{b-i}$. Similarly,
 we have  $|N(B)|=\sum_{i=0}^{t-1}\binom{b}{i}\binom{n-b}{a-i}$.

 It is well known that for each $A\in\binom{[n]}{k}$, the
stabilizer of $A$ is a maximal subgroup of $S_n$ subject to $n\neq
2k$ \cite{sym}. Therefore, the action of $S_n$ on $\binom{[n]}{k}$
is imprimitive if and only if $n=2k\geq 4$, and the only imprimitive
sets are all  pairs of complementary subsets. If $n=2a\geq 4$ and
$(n,t)\neq (a+b,1)$, then from $a+b<n+t$ it follows $b<a+t$. For
every pair $A$ and $\overline {A}$ in ${[n]\choose a}$, it is easy
to verify that $\{C\cup\{i\}:C\in {\overline  A\choose a-1}, i\in
A\}\subseteq N(A)\backslash N(\overline {A})$ if $b=a$ and $t>1$,
${\overline  A\choose b}\subseteq N(A)\backslash N(\overline {A})$
if $b<a$, and $\{\overline  A\cup C:C\in { A\choose b-a}\}\subseteq
N(A)\backslash N(\overline {A})$ if $b>a$. This implies
$|N(A)\backslash N(\overline {A})|>1$ hence $|N(A)\cap N(\overline
{A})|<|N(A)|-1$. Therefore, $\{A,\overline  A\}$ is not a fragment
for every $A\in {[n]\choose a}$. We thus prove that every fragment
in $\mathcal X$ and $\mathcal Y$ is primitive. Then, by Theorem
\ref{Gd}, inequality (\ref{cross1}) holds.

To complete the proof of Theorem \ref{subsets} we need to determine
all nontrivial fragments.  Suppose there is a nontrivial fragment
 in $\binom{[n]}{a}$ or $\binom{[n]}{b}$. Without
loss of generality we assume that $\mathcal S$ is a minimal-sized
one in $\binom{[n]}{a}$. By Theorem \ref{Gd},
$\binom{n}{a}=\binom{n}{b}$, i.e., $b=a$ or $b=n-a$. Clearly, $S_n$
is not isomorphic to a subgroup of $D_{n!}$ for $n\geq 4$.
Therefore, by Proposition \ref{2frag}, there are no 2-fragment in
$\mathcal{F}(\mathcal{X})$ and $\mathcal{F}(\mathcal{Y})$, which
implies that $\mathcal S$ is balanced.

For each $C\subseteq [n]$, $S_C$ is embedded into $S_n$ in a natural
way:  for $\sigma\in S_C$, let $\sigma$ fixes elements of $\overline
C$. Now, take a $C\in \mathcal S$ and let $\Gamma=S_C\times
S_{\overline C}$ and $\Gamma_{\mathcal S}=\{\sigma\in \Gamma:
\sigma(\mathcal S)=\mathcal S\}$. Then $C\in\sigma(\mathcal S)$ for
each $\sigma\in\Gamma$. Since $\mathcal S$ has more than one
elements, we have $\Gamma_{\mathcal S}\neq \Gamma$. Otherwise,
$S_C\times S_{\overline C}$ and $S_B\times S_{\overline B}$ (for
some $B\in\mathcal S\backslash\{C\}$) will generate whole $S_n$ so
that $ S={[n]\choose a}$, yielding a contradiction. Then, by
Proposition \ref{3} we have that $[\Gamma:\Gamma_{\mathcal S}]$, the
index of $\Gamma_{\mathcal S}$ in $\Gamma$, equals 2. Now, let
$\Gamma_{\mathcal S}[C]$ be the projection of $\Gamma_{\mathcal S}$
onto  $S_C$. Then, $\Gamma_{\mathcal S}[C]$ is a subgroup of $S_C$
of index $\leq 2$. That is, $\Gamma_{\mathcal S}[C]=S_C$ or $A_C$.
From this we see that $A_C\times S_{\overline C}$ and $S_C\times
A_{\overline C}$ are the only index-2 subgroups of $\Gamma$. That
is,  $\Gamma_{\mathcal S}=A_C\times S_{\overline C}$ or  $S_C\times
A_{\overline C}$. For any $B\in \mathcal S\backslash\{C\}$,
$a=|B\cap C|+|B\cap\overline C|$.
If  $|B\cap C|>1$, let $(i,j)$ be an interchange, where $i,j\in
B\cap C$. Then, $(i,j)$ fixes both $C$ and $B$. The
semi-imprimitivity of $\mathcal S$ implies $(i,j)\in\Gamma_{\mathcal
S}$. This yields  $\Gamma_{\mathcal S}=S_C\times A_{\overline C}$.
From this process it follows that, for each $B\in\mathcal S$, there
exists at most one of $|B\cap C|$ and $|B\cap\overline C|$ being
grater than 1. Note that if $B\subseteq \overline C$, then $S_C$ and
$S_B$ fix both $C$ and $B$, i.e., $S_C\times S_B\subseteq
\Gamma_{\mathcal S}$.  It is clear, however, that neither $A_C\times
S_{\overline C}$ nor $S_C\times A_{\overline C}$ contain $S_C\times
S_B$. We therefore obtain that $|C\cap B|=1$ for every $B\in\mathcal
S$, or $|C\cap B|=a-1$ for every $B\in\mathcal S$.

Suppose $|C\cap B|=1$ for every $B\in \mathcal S$. Without loss of
generality we assume $C\cap B=\{1\}$ for some $B\in \mathcal S$. In
this case, if $a>2$, then $|B\cap\overline C|\geq 2$, so
$\Gamma_{\mathcal S}=A_C\times S_{\overline C}$. On the other hand,
we can find distinct $i,j\in C$ such that
$(1,i,j)(B)=B\backslash\{1\}\cup\{i\}\in \mathcal S$ because
$(1,i,j)\in A_C$. From this it follows that  $(1,i)(\mathcal S)$
contains  more than one element of $\mathcal S$, hence
$(1,i)\in\Gamma_{\mathcal S}$. The contradiction proves $a=2$. Thus
$\mathcal S$ consists of all 2-subsets $\{1,i\}$'s for $i\in [2,n]$.
Since $t<\min\{a,b\}$ and $(n,t)\neq (a+b,1)$, we have $t=1$ and
$b=2$. Then  $d(\mathcal X)={n-2\choose 2}$ and $N(\mathcal
S)={[2,n]\choose 2}$ satisfying $|N(\mathcal S)|-|\mathcal
S|=d(\mathcal X)-1$, that is, $\mathcal S$ is a fragment in
${[n]\choose2}$.

Suppose now $|C\cap B|=a-1>1$ for every $B\in \mathcal S$.
In this case, we may similarly prove that $n-a=2$, $b=a$, $t=n-3$
and $\Gamma_{\mathcal S}=S_{C}$. Thus $\mathcal S=\{\sigma(B\cap
C)\cup\{i\}: \sigma\in S_C\ \mbox{ and $\{i\}=B\cap\overline
C$}\}\cup\{C\}={A\choose n-2}$ where $A=B\cup C$ is a $(a+1)$-subset of
$[n]$. It is easy to verify that $\mathcal S$ is a fragment in
${[n]\choose n-2}$. \qed

\section{Proof of Theorem \ref{subspaces}}
Similarly to the proof of Theorem \ref{subsets}, put $\mathcal
{X}=\qchoose{V}{a}$ and $\mathcal {Y}=\qchoose{V}{b}$. The bipartite
graph $G(\mathcal{X},\mathcal{Y})$ is defined by the
cross-$t$-intersecting relation between  $\mathcal{X}$ and
$\mathcal{Y}$:  for $A\in\mathcal{X}$ and $B\in\mathcal{Y}$,
$(A,B)\in E(G)$ if and only if $\dim(A\cap B)<t$. Analogously to the
families of sets, $G(\mathcal{X},\mathcal{Y})$ is connected and
 non-complete.
 Let $GL(V)$ denote the general linear group of
$V$, which consists of all invertible linear transformations of $V$.
Clearly,
 $GL(V)$ transitively acts on $\mathcal {X}$ and $\mathcal {Y}$, respectively, in a natural way, and preserves the
 cross-$t$-intersecting relation. So $d(\mathcal {X})=|N(A)|$ for
 $A\in \mathcal {X}$. It is easy to see that
 \[N(A)=\bigcup_{i=0}^{t-1}N_i^b(A),\]
 where $N_i^b(A)=\left\{T\in \qchoose{V}{b}:\dim(T\cap A)=i\right\}$.
 To determine $|N(A)|$, we need a usefull result, stated as a lemma as follows.
\begin{lemma}\label{chen-rota} (\cite[Proposition 2.2]{chen}) Let $A\in \qchoose{V}{a}$.
 Then $|N_0^b(A)|=\qchoose{n-a}{b}q^{ab}$.
\end{lemma}

Suppose $i>0$ and let $C$ be an arbitrary $i$-subspace of $A$. We
consider the quotient space $V/C$. Then $\dim(V/C)=n-i$,
$\dim(A/C)=a-i$ and
$|N_0^{b-i}(A/C)|=\qchoose{n-a}{b-i}q^{(a-i)(b-i)}$. So
$|N_i^b(A)|=\qchoose{a}{i}|N_0^{b-i}(A/C)|=\qchoose{a}{i}\qchoose{n-a}{b-i}q^{(a-i)(b-i)}$
(See also \cite[Lemma 4]{wy} and \cite[Lemma 2.4]{wzhjt}).
 Thus
 $$|N(A)|=\sum_{i=0}^{t-1}q^{(a-i)(b-i)}\qchoose{a}{i}\qchoose{n-a}{b-i}=d(\mathcal
 X).$$  Similarly, $d(\mathcal
 Y)=\sum_{i=0}^{t-1}q^{(a-i)(b-i)}\qchoose{b}{i}\qchoose{n-b}{a-i}$.

For a subspace $A$ of $V$, by $GL(V|A)$ we denote the stabilizer of
$A$ in $GL(V)$. It is well known that for $A\in \mathcal {X}$,
$GL(V|A)$ is a maximal subgroup of $GL(V)$ \cite{Asch}, so the
action of $GL(V)$ on $\mathcal {X}$ is primitive.  Then, by Theorem
\ref{Gd}, inequality (\ref{cross2}) holds, and each nontrivial
fragment is a semi-imprimitive set under the action of $GL(V)$.

To complete the proof of Theorem \ref{subsets} we need to determine
all nontrivial fragments.  Suppose there is a nontrivial fragment
 in $\mathcal X$ or $\mathcal Y$. Without
loss of generality we assume that $\mathcal S$ is a minimal-sized
one in  $\mathcal X$.  By Theorem \ref{Gd},
$\qchoose{n}{a}=\qchoose{n}{b}$, i.e., $b=a$ or $b=n-a$. Clearly,
$GL(V)/K$ is not isomorphic to a subgroup of $D_{|\mathcal X|}$,
where $K$ is the kernel of the action of $GL(V)$ on $\qchoose{V}{a}$
or $\qchoose{V}{b}$. Therefore, by Proposition \ref{2frag}, there
are no 2-fragment in $\mathcal{F}(\mathcal{X},\mathcal{Y})$, which
implies that $\mathcal S$  is balanced.

Take  $C\in \mathcal S$, write $\Gamma=GL(V|C)$ and
$\Gamma_{\mathcal S}=\{\sigma\in\Gamma:\sigma(\mathcal S)=\mathcal
S\}$. Then, $\Gamma\neq\Gamma_{\mathcal S}$, and  again by
Proposition \ref{2frag}, $[\Gamma, \Gamma_{\mathcal S}]=2$ so that
$\Gamma=\Gamma_{\mathcal S}\cup\gamma\Gamma_{\mathcal S}$ for some
$\gamma\in\Gamma$. Thus $\mathcal S$ and $\gamma(\mathcal S)$ are
the only nontrivial fragments containing $C$. From the structure of
$\Gamma$ it follows that $\Gamma$ is transitive on  $N_i^a(C)$ for
each $i=0,1,\ldots, a$, whenever   $N_i^a(C)\neq\emptyset$.
 Set $\mathcal
S_i=\mathcal S\cap N_i^a(C)$ for  $0\leq i<a$. If $\mathcal
S_i\neq\emptyset$, then
\[N_i^a(C)=\mathcal S_i\cup\gamma(\mathcal S_i)=\left\{L+R:\mbox{ $L\in
\qchoose{C}{i}$ and $R\in N_{0}^{a-i}(C)$}\right\}.\]  From this we
see that  $\Gamma_{\mathcal S}$ is transitive on $\mathcal S_i$. It
is clear that the restriction of $\Gamma$ on $C$ is $GL(C)$.
Therefore, the induced action of $\Gamma$ on $\qchoose{C}{i}$ is
primitive, thus the action of $\Gamma_{\mathcal S}$ on
$\qchoose{C}{i}$ is transitive. This means that if $L_0+R_0\in
\mathcal S_i$ for some $R_0\in N_{0}^{a-i}(C)$, then  $L+R_0\in
\mathcal S_i$ for every $L\in \qchoose{C}{i}$. We complete the proof
by two cases.

Case 1: $n-a=a-i$. Suppose that $L_0+R_0\in \mathcal S_i$ and
$\{\alpha_1,\ldots,\alpha_{a-i}\}$ is a basis of $R_0$. Then bases
of elements of $N_{0}^{a-i}(C)$ are of the form
$\{\alpha_1+\beta_1,\ldots,\alpha_{a-i}+\beta_{a-i}\}$, where
$\beta_i\in C$. Put $Q=\{R\in N_{0}^{a-i}(C): L+R\in \mathcal S_i\
\mbox{for some}\ L\in  \qchoose{C}{i}\}$. Then $R_0\in Q$. For given
$\beta_1,\ldots,\beta_{a-i}\in C$, let $R_j$ be the subspace
generated by $\alpha_1+\beta_1,\ldots,\alpha_{j}+\beta_{j},
\alpha_{j+1},\ldots,\alpha_{a-i}$. Assume $R_j\in Q$. Then the above
discussion implies $L+R_j\in Q$ for every $L\in \qchoose{C}{i}$.
Thus, we can take an $L\in\qchoose{C}{i}$ containing $\beta_{j+1}$
so that $L+R_j=L+R_{j+1}$, that is, $R_{j+1}\in Q$.  This proves
$\mathcal S_i=N_i^{a}(C)$, yielding a contradiction.

Case 2: $n-a>a-i$. Consider the natural map $\nu$ from $V$ onto the
quotient space $V/C$, that is, $\nu(A)=(A+C)/C$, written as $\bar
A$, for any subspace $A$ of $V$. Then $\nu(N_i^{a-i}(C))=\qchoose
{V/C}{a-i}$. It is clear that $\Gamma$ acts on $V/C$ and  $\Gamma/K$
is isomorphic to $GL(V/C)$, where $K$ is the kernel of the action.
Then the primitivity of the action implies that $\Gamma_{\mathcal
S}K/K$ is transitive on $\qchoose {V/C}{a-i}$. This means that for
each $\bar R_0 \in \qchoose {V/C}{a-i}$, there is an $R_0\in
N_{0}^{a-i}(C)$ such that $\nu(R_0)=\bar R_0$ and $L_0+R_0\in
\mathcal S_i$ for some $L_0\in \qchoose{C}{i}$. Then, by Case 1 we
prove $\mathcal S_i=N_i^{a}(C)$, yielding a contradiction, again.

We thus prove that the graph has no nontrivial fragments.\qed

\section{Proof of Theorem \ref{group}}
We first prove a general result. Let $\Gamma$ be a transitive
permutation group on $\Omega$ with the identity 1. By the group and
a positive integer $t$ with $1\leq t\leq |\Omega|-2$ we define a
simple graph, written as $G_t=G_t(\Gamma)$, whose vertex set is
$\Gamma$, and whose edge set consists of all pairs $\sigma\tau$ such
that $|\{x\in \Omega:\sigma(x)=\tau(x)\}|< t$. Let $\Gamma_L$ and
$\Gamma_R$ denote the left and right regular action on $\Gamma$,
respectively. Then $\Gamma_L\times \Gamma_R$ (not necessarily a
direct product) induces an automorphism group of $G_t(\Gamma)$. In a
natural way, we can view $G_t(\Gamma)$ as a bipartite graph
$G_t(\Gamma,\Gamma)$, which is  part-transitive under the action of
$\Gamma_L$ and $\Gamma_R$.

\begin{lemma}\label{imp} Suppose that $A$ is an imprimitive set in $\Gamma$
under the action of  $\Gamma_L\times \Gamma_R$. Then $A$ is a coset
of a non-trivial normal subgroup of $\Gamma$.
\end{lemma}
\textbf{Proof.} Since $A$ is an imprimitive set, we have that
$1<|A|<|\Gamma|$, and for every $\alpha\in \Gamma$, $\alpha A$ is
also an imprimitive set. Without loss of generality we assume that
$1\in A$. From this it follows that $\alpha\in \alpha A$ and $1\in
\alpha^{-1}A$ for
 each $\alpha\in A$, hence $\alpha A=\alpha^{-1}A=A$, which implies
 that $A$ is a subgroup $\Gamma$. Furthermore, for every
 $\gamma\in\Gamma$, $1\in (\gamma^{-1}A\gamma)\cap A$, hence  $\gamma^{-1}A\gamma=
 A$, proving that $A$ is a normal subgroup of $\Gamma$. \qed

We now consider the  graph $G_t(S_n)$ where $n\geq 4$ and $1\leq
t\leq n-2$. For $0\leq i<n$, by $\mathcal D_n^i$ we denote the set
of all permutations in $S_n$ which have exact $i$  fixed points. The
elements of  $\mathcal D_n^0$ are known for the derangements of
$[n]$. As usual, set $|\mathcal D_n^0|=D_n$. By definition,
$G_t(S_n)$ is the Cayley graph on $S_n$ generated by $ \mathcal
G_t$, where  $\mathcal G_t=\cup_{i=0}^{t-1}\mathcal D_n^i$.  (cf.
\cite{ku-per}). It is not difficult to compute that for every
$\sigma\in S_n$,
\[|N(\sigma)|=|\mathcal G_t\sigma|=|\mathcal G_t|=\sum_{i=0}^{t-1}{n\choose
i}D_{n-i}.\]

Let $\mathcal{S}$ be a fragment in $S_n$. Then for any $\sigma\in
S_n$, $\sigma\mathcal{S}$ is also a fragment. Without loss of
generality, we assume that $1\in \mathcal{S}$ and set
$\mathcal{S}^*=\mathcal{S}\backslash\{1\}$. By definition we have
that $|N(\mathcal{S})|=|\mathcal G_t\mathcal{S}|=|\mathcal
G_t|+|\mathcal{S}^*|$, that is
\begin{equation}\label{ll}
|\mathcal G_t\mathcal{S}^*\backslash\mathcal G_t|=|\mathcal{S}^*|.
\end{equation}

If $\mathcal S$ is imprimitive, then Lemma \ref{imp} implies that
$\mathcal S$ is a nontrivial normal subgroup of $S_n$. It is well
known that the only nontrivial normal subgroups of $S_n$ are $A_n$
and the quaternary group $V_4=\{1, (12)(23),(13)(24),(14)(23)\}$ for
$n=4$. Since $A_n$ has  index 2 in $S_n$ and $\mathcal
G_t\not\subset A_n$, $\mathcal G_t A_n=S_n$, hence $A_n$ is not a
fragment in $S_n$ for $n\geq 4$. And, for $n=4$ and $t=1,2$, it is
straightforward to verify that $|\mathcal G_tV_4^*\setminus\mathcal
G_t|>3$,  so $V_4$ is not a fragment in $S_4$. We thus prove that
every fragment in $S_n$ is primitive. Then, by Theorem \ref{Gd} we
obtain inequality (\ref{cross3}). Moreover,  by Proposition
\ref{2frag}, it is easy to verify that $G_t(\Gamma,\Gamma)$ has no
2-fragments.

Suppose that there is a nontrivial fragment $\mathcal S$ in $S_n$.
Then, by Proposition \ref{3},  $\mathcal S$  is balanced and
$|\mathcal S|>2$. Without loss of generality we may assume
$1\in\mathcal S$.
 Set $H=\{h\in S_n:
h\mathcal{S}=\mathcal{S}\}$. Clearly, $H$ is a subgroup of $S_n$. If
$H=\{1\}$, then $\sigma\mathcal{S}=\tau\mathcal{S}$ implies
$\sigma=\tau$ for any $\sigma,\tau\in S_n$, hence for any distinct
$a,b\in\mathcal{S}$,
 by the semi-imprimitivity of $\mathcal{S}$, we have $a^{-1}\mathcal{S}\cap b^{-1}\mathcal{S}=\{1\}$.
  We thus obtain more than 2
$|\mathcal S|$-fragments  containing 1, contradicting
 Proposition \ref{3}.  Therefore, $|H|>1$ and $S=\cup_{b\in S}Hb$. For each $a\in\mathcal{S}$,
  it is evident that $Ha\subset \mathcal{S}\cap \mathcal{S}a$. So the semi-imprimitivity of $\mathcal{S}$ implies
  that $\mathcal{S}=\mathcal{S}a$, which implies that
  $\mathcal{S}$ is a subgroup of $S_n$. We have seen that $\mathcal S$ is not
  normal. i.e., there is a $\sigma\in S_n$ with $\sigma^{-1}\mathcal S\sigma\neq\mathcal
  S$. However, each $\sigma^{-1}\mathcal S\sigma$ contains 1.   Again by Proposition \ref{3},  the
normalizer $N_{S_n}(\mathcal S)$ is an index-2 subgroup of $S_n$,
i.e., $N_{S_n}(\mathcal S)=A_n$ because $A_n$ is the only index-2
subgroup of $S_n$. So $\mathcal S$ is a normal subgroup of $A_n$. It
is well known that $A_n$ is a simple group for $n\geq 5$, therefore
$A_n$ has no nontrivial normal subgroup for $n\geq 5$, and $A_4$ has
the only nontrivial normal subgroup $V_4$. We have seen that neither
 $A_n$ nor $V_4$ are fragments of $G_t(\Gamma,\Gamma)$. We thus prove that the graph has no nontrivial fragments.
This completes the proof. \qed

\end{document}